\documentclass[12pt]{article}
 \usepackage{verbatim}

\date{}
\title{Regularity for  radial solutions of degenerate fully nonlinear equations.  }

\author{ I. Birindelli, F. Demengel}

\newtheorem{theo}{Theorem}[section]
\newtheorem{prop}[theo]{Proposition}
\newtheorem{rema}[theo]{Remark}
\newtheorem{defi}[theo]{Definition}
\newtheorem{cor}[theo]{Corollary}
\newtheorem{lemme}[theo]{Lemma}

\newcommand{\ph}{\varphi}
\newcommand{\R}{{\rm I}\!{\rm  R}}
\def\R{{\rm I}\!{\rm  R}}

\def\la1{\lambda_1}
\def\grad{\nabla}
\def\ph1{\varphi_1}
\newcommand{\N}{{\bf N}}
\newcommand{\Z}{{\bf Z}}

\begin{document}
\maketitle
\begin{abstract}
In this paper we prove the ${\cal C}^{1, \beta}$ regularity of the solutions of radial solutions for fully nonlinear degenerate equations.
\end{abstract}

\section{Introduction}

In this paper we prove the regularity of radial solutions of 
$$ F(x,\grad u,D^2 u)=f,$$
where $F$ is  a fully nonlinear degenerate elliptic operator, homogenous of degree 1 in the Hessian and homogenous of some degree $\alpha>-1$ in the gradient, 
which is elliptic when the gradient is not null.
Precise conditions on $F$ will be stated in the next section.

 The class of operators we consider includes:
$$F(\grad u,D^2 u)= |\nabla u|^\alpha {\cal  M}_{a,A}(D^2u)$$ 
 where $ {\cal  M}_{a,A}$ is one of the Pucci operators (i.e. either ${\cal  M}_{a,A}={\cal  M}_{a,A}^+$ or ${\cal  M}_{a,A}={\cal  M}_{a,A}^-$) , 
 $$F(\grad u,D^2 u)=\Delta_{\alpha+2}u=\rm{div}( |\nabla u|^\alpha\nabla u)$$
 or, more in general,
$$F(\grad u,D^2 u)= |\nabla u|^\alpha( p_1tr(D^2u)+ p_2 \langle D^2u {\grad u\over |\grad u|}, {\grad u\over |\grad u|} \rangle )$$
with $p_1>0$ and $p_1+p_2>0$.

In a previous paper  \cite{BD2} we proved that for $\alpha\in (-1,0]$ all solutions  are ${\cal C}^{1,\beta}$ if $F$ satisfies

\begin{equation}\label{00} F(\grad u,D^2 u)= |\nabla u|^\alpha\tilde F(D^2 u),
\end{equation}
( if  $\tilde F$ is concave we obtained that the solutions are ${\cal C}^{2,\beta}$).

Here we prove that even for $\alpha>0$, the radial solutions are ${\cal C}^{1,\frac{1}{1+\alpha}}$ 
everywhere, and if the dependence on the Hessian is convex, in points where the derivative is not zero, the solutions are ${\cal C}^{2,\beta}$. Observe that where the radial derivative is zero the H\"older continuity of the first derivative is optimal. Indeed it is easy to see that
$ u(r) = r^{\alpha+2\over \alpha+1}$
is a viscosity solution of 
 $$ |\nabla u|^\alpha {\cal M}_{a,A}^+ (   D^2 u) =c$$
 for $c=\left( {\alpha+2\over \alpha+1}\right)^{\alpha+1}  A( {1\over 1+\alpha} + N-1). $

Beside its intrinsic interest, the regularity question was raised naturally while proving the simplicity of the principal eigenfunctions.
In recent years, the concept of principal eigenvalue has been extended to fully nonlinear operators, by means of the maximum principle (see \cite{BNV}).
The values
$$\lambda^+(\Omega)=\sup\{\lambda, \exists \phi>0 \ \mbox{in}\ \Omega,  F(x,\grad\phi, D^2 \phi)+\lambda \phi^{1+\alpha} \leq 0\quad\mbox{in}\ \Omega\}$$
$$\lambda^-(\Omega)=\sup\{\lambda, \exists \psi<0 \ \mbox{in}\ \Omega, F(x,\grad\psi, D^2 \psi)+\lambda |\psi|^{\alpha}\psi \geq 0\quad\mbox{in}\ \Omega\}$$
are  generalized eigenvalues in the sense that there exists a non trivial solution to the Dirichlet problem 
$$ F(\grad\phi, D^2 \phi)+\lambda^{\pm}(\Omega) |\phi|^{\alpha}\phi= 0\quad\mbox{in }\quad \Omega,\ \phi=0 \quad\mbox{on}\quad \partial\Omega.$$
In \cite{BD2}, we proved that, for $F$ satisfying (\ref{00}), these eigenfunctions are simple as long as $\partial\Omega$ has only one connected component. This result   extends to the situation where $\partial \Omega$ has  at most two connected components  when the dimension is  $2$, the proof uses the fact that the eigenfunctions are ${\cal C}^{1,\beta}$.

Let us emphasize that regularity results for degenerate elliptic operators that are not in divergence form are in general very difficult. The difficulty  comes from the fact that difference of solutions are not sub or super solutions of some elliptic equation. As an example,  let us recall that for the  infinity Laplacian $\Delta_{\infty}$ the solutions are known to be in
 ${\cal C}^{1,\beta}$ for small $\beta>0$ only in dimension 2, and only for $f\equiv 0$, see \cite{ES}. 
 
 On the other hand, N. Nadirashvili, S. Vladut in  \cite{NV} prove $C^2$ regularity of radial solutions for a large class of fully nonlinear operators uniformly elliptic.

 \bigskip
 Of course we use viscosity solutions, and since, as can easily be imagined, the difficulties arise where the derivative is zero, our first concern is to check that if $u^\prime\neq 0$ in the viscosity sense at some point,  then this holds in a neighborhood, and furthermore, in that neighborhood the solution is   ${\cal C}^{1}$. Then we treat the points where 
 $u^\prime= 0$ in the viscosity sense.
The proof relies  only on the comparison principle,  Hopf principle,   the regularity results of \cite{CDV} and \cite{CaC2} , together with some classical analysis.

 \section{Hypothesis and known results.}
 In all the paper we suppose that $\Omega$ is a  ball or an annulus.
We shall consider solutions of the following equation
\begin{equation}\label{Theq}
F(x, \grad u,D^2u)=f(|x|).
\end{equation}

The operator $F$ is supposed to  be continuous on $\Omega \times (\R^N)^\star \times S$,  where $S$ is the space of $N\times N$ symmetric matrices  and to satisfy:
\begin{enumerate}
\item[(H1)] For some $\alpha >-1$, for all $x\in \Omega$,  for all $p\neq 0$ and $N\in S$ and for all $t\in \R$ and $ \mu>0$,
$F(x, tp, \mu N) = |t|^\alpha \mu F(x, p, N)$.

\item[(H2)] $F$ is  fully nonlinear elliptic, i.e  there exist some positive constants $a$ and $A$,  such that for any $M\in S$ and $N\geq 0$ in $S$, one has 
  $$a|p|^\alpha  tr(N)\leq  F(x,p,  M+N)- F(x, p,M) \leq A|p|^\alpha  tr(N).$$
  
\item[(H3)] Furthermore $F$ is an Hessian operator, i.e. for any $M\in S$ and $O\in {\cal O} (n, \R)$, $ F(Ox, ^t Op, ^t O M O) =  F(x, p, M)$. 

\item [(H4)]
There exists $\nu>0$ and $\kappa \in ]1/2,1]$ such that for all $|p|=1$ ,
$|q|\leq {1\over 2}$, $M\in {\mathcal S}$
 $$|F(x,p+q,M)- F(x,p,M)|\leq \nu |q|^\kappa |M|$$
which implies by homogeneity 
that for all  $p\neq 0$ ,
$|q|\leq {|p|\over 2}$, $M\in {\mathcal S}$
 $$|F(x,p+q,M)-F(x,p,M)|\leq \nu |q|^\kappa|p|^{\alpha-\kappa} |M|$$

\end{enumerate}
We need to precise what we mean by viscosity solutions : 
       \begin{defi}\label{defvs}
Let $\Omega$ be a domain in
$\R^N$, let $g$ be a continuous function on $\Omega\times \R$, then
$v$,   continuous  in $\Omega $  is called a viscosity super-solution (respectively sub-solution) 
of
$F(x, \grad u,D^2u)=g(x,u)$ if for all $x_0\in \Omega$,

- Either there exists an open ball $B(x_0,\delta)$, $\delta>0$  in $\Omega$
on which 
$v$ equals to   a constant $c$ and 
$0\leq g(x,c)$, for all $x\in B(x_0,\delta)$
(respectively 
$0\geq g(x, c)$ for all $x\in B(x_0,\delta)$).

- Or
 $\forall \varphi\in {\mathcal C}^2(\Omega)$, such that
$v-\varphi$ has a local minimum (respectively local maximum) at  $x_0$ and $\grad\varphi(x_0)\neq
0$, one has
$
F(x_0,  \grad\varphi(x_0),
 D^2\varphi(x_0))\leq g(x_0,v(x_0)).
$
(respectively  
$F(x_0,  \grad\varphi(x_0),
 D^2\varphi(x_0))\geq g(x_0,v(x_0))).$

A	 viscosity solution is a function which is both a super-solution and a sub-solution.
\end{defi}
 Let us observe that in the case where $\alpha >0$, the operator is well defined everywhere, and then it is a natural question to ask if the viscosity solutions in the sense above are the same as the viscosity solutions in the usual sense. The answer is yes as is proved in the appendix of this paper.

 From now on we suppose that  $\alpha \geq 0$ and  that the solutions are "radial ".  Let us observe that the hypothesis that $F$ be a fully nonlinear elliptic 
  hessian operator  implies  that 
 there exists  some  operator $H$ defined on $\R^+\times \R^2$, such that   if $u(x) = g(r)$ is radial and ${\cal C}^2$,  
 $F(x,\nabla u, D^2 u) =  H(r, g^{\prime\prime},g^\prime)$
  with 
  \begin{eqnarray*}
  |g^\prime|^\alpha \left(A \left((g^{\prime\prime})^- + {N-1\over r} (g^\prime) ^-\right)\right.&+& \left. a  \left((g^{\prime\prime})^+ + {N-1\over r} (g^\prime) ^+\right) \right)\\
  &\leq&  H(r,g^{\prime\prime}, g^\prime)\\
  &\leq& 
  |g^\prime|^\alpha \left(A \left((g^{\prime\prime})^+ + {N-1\over r} (g^\prime )^+\right)\right.\\
  &+ & a\left. \left( (g^{\prime\prime})^- + {N-1\over r} (g^\prime) ^-\right)\right)
  \end{eqnarray*}

   \bigskip

\noindent We now recall some known results  concerning the  operators  considered : 

\begin{prop}[\cite{BD1}]\label{comp1}
Suppose that $\Omega$ is a bounded  open set. 
 Suppose that $f$ and $g$ are continuous on $\overline{\Omega}$ and $f\geq  g$. Assume  that $\beta$ is some continuous and non decreasing function  on $\R$, and that  $u$ and $v$ are continuous   respectively sub- and super-solutions   of the equation 
$$F(x,\nabla u, D^2 u) -\beta (u) \geq  f$$
$$F(x, \nabla v, D^2 v)-\beta (v) \leq  g$$
 with $u\leq v$ on $\partial \Omega$. Then if $ f> g$ in $\Omega$,  or if $f\geq g$ but $\beta$ is increasing, 
  $u\leq v$ inside $\Omega$.
  \end{prop}
  
We shall also need the Lipschitz regularity of the solutions: 
   
    \begin{prop}[\cite{BD1}]\label{lip}
    Suppose that $\Omega$ is an open  bounded regular domain of $\R^N$. 
  Suppose that $u$ is a function in ${\cal C}(\overline{\Omega})$ which is a solution of 
  $$\left\{ \begin{array}{lc}
  F(x, \nabla u, D^2 u) = f&\ {\rm in}\  \Omega,\\
  u=0 & {\rm on } \ \partial \Omega,
  \end{array}\right.$$
 then $u$ is Lipschitz continuous.
  \end{prop}
  
  \begin{rema}
   We shall use this Proposition for radial solutions, so we shall fix   $r_o$, and use the previous Proposition  for  $u-u(r_o)$ which is a solution   of the above homogeneous Dirichlet problem on $B(0, r_o)$.
\end{rema}

In all the paper, $f$ denotes some continuous and radial function. 
\section{${\cal C}^1$ regularity}
 In all the sequel we denote for simplicity 
 $$F[u]:= H(r, u^{\prime\prime},\frac{u^\prime}{r})$$ and  $u$ is supposed to be a continuous  radial solution of $F[u]=f(r)$ on either the interval $[0, R)$ or the interval $(R_1,R_2)$. 
\begin{defi}
For any $(p,q)\in \R^2$, we define the paraboloid 
$$w(p, q,  r)(s)= u(r) + p(s-r) + {q\over 2} (s-r)^2.$$
\end{defi}
We also give the following
\begin{defi}
For a Lipschitz continuous function $u$, we define the following  so called   derivative numbers of $u$ : 
 $$\lambda_g(r_1) = \liminf _{r\rightarrow r_1, r< r_1}{u(r)-u(r_1)\over r-r_1} , $$     
 $$\Lambda_g(r_1) = \limsup _{r\rightarrow r_1, r< r_1}{u(r)-u(r_1)\over r-r_1},$$
$$\lambda_d(r_1) = \liminf _{r\rightarrow r_1, r> r_1}{u(r)-u(r_1)\over r-r_1},  $$
  $$\Lambda_d(r_1) = \limsup _{r\rightarrow r_1, r> r_1}{u(r)-u(r_1)\over r-r_1} .$$ 
Clearly for $r_1=0$, only the "right" derivatives are defined.

We shall say that $u^\prime (\bar r)>0$ in the viscosity sense   (respectively $u^\prime(\bar r) <0$) if $\inf(\lambda_g(\bar r), \lambda_d(\bar r))>0$  (respectively if 
$\sup (\Lambda_d(\bar r), \Lambda_g(\bar r) )<0$). 

 On the opposite we shall  say that $u^\prime(\bar r)=0$ in the viscosity sense if one has 
 $$\lambda_g (\bar r)\lambda_d (\bar r) \leq 0\,  \ \mbox{and}\,\   \Lambda_d (\bar r)\Lambda_g(\bar r) \geq 0.$$ 
\end{defi}
\begin{rema}

 Let us note first that    all the numbers  defined above exist and are finite for $u$ a solution of (\ref{Theq}) since the solutions are known to be Lipschitz. Furthermore we proved in \cite{BDr} that since $f$ is bounded,  
$\Lambda_g \geq  \lambda_d$, $\Lambda_d \geq  \lambda_g$ .  Finally, if  all these numbers coincide on $\bar r$, $u^\prime$ exists  on $\bar r$ in the classical sense. 
\end{rema}

We begin with a simple lemma

\begin{lemme}\label{lemuprime0}
Suppose that $u$ is a  radial continuous viscosity solution of 
$F[u] = f$  on $[0, R)$,  then $u^\prime (0)$ exists and it is zero.
\end{lemme}

 \begin{rema}
  
   We want to point out that for radial function, i.e. for the continuous functions $u$ defined on  some ball of $\R^N$, such that there exists $v$ continuous on $[0, r]$, with 
$u(x) = v(|x|)$, in order to test  on points $x\neq 0$, it is sufficient to use test functions 
which are radial. A consequence of Lemma \ref{lemuprime0} and Definition \ref{defvs} is that we do not need to test at the point zero.  As a consequence $u$ is a viscosity  supersolution  of $F[u] = f$ 
in $B(0, R)$ if and only if   $u^\prime(0)$ exists and is zero, and for all $r\neq 0$
and for all  $\varphi$ which is  ${\cal C}^2$ around $\bar r \neq 0$ which touches   $u$ by below  on $\bar r$
$$F[\varphi](\bar r) \leq f(\bar r).$$
\end{rema}
 {\em Proof of Lemma \ref{lemuprime0}.} We want to prove that 
 $$\Lambda_d(0)=\lambda_d(0)= 0.$$ 
Suppose that 
$\Lambda_d(0) = m >0$.  Let $m^\prime< m$ arbitrary close to it.  Choose $\delta$ small enough in order that 
$$ a (m^\prime)^{1+\alpha} {N-1\over \delta}> |f|_\infty.$$  
By hypothesis, for such $\delta$,  there exists $r\in ]0, \delta]$ such that 
$${u(r)-u(0)\over r} \geq m^\prime \ \mbox{i.e.}\ u(r)\geq u(0)+m^\prime r.$$ 
Let  $w:= w(m^\prime, 0, 0)$ so that $w(0) = u(0)$, $w(r) \leq u(r)$,  and   
$$F[w] (s)\geq a (m^\prime)^{1+\alpha } {N-1\over s} \geq f(s)\ \mbox{on}\ [0, \delta].$$
Then, by Proposition \ref{comp1}, $w(s) \leq u(s)$,  in $ [0, r]$. 
Hence   $${u(s)-u(0)\over s} \geq m^\prime,\  \mbox{and}\ \lambda_d(0)  \geq m^\prime.$$ 
This implies that $\Lambda_d(0)=\lambda_d(0)$.
        
        We now suppose that $\Lambda_d(0) \leq 0$, so either  $\lambda_d(0)  = 0$ and then they are equal and we are done, or $\lambda_d(0) =- m <0$.  

Let $0<m^\prime< m$ arbitrary close to it. Choose $\delta$ as above, 
by hypothesis, there exists $r\in ]0, \delta]$,  such that 
$${u(r)-u(0)\over r} \leq -m^\prime .$$ 

Let  $w\equiv w(- m^\prime, 0, 0)$, then $w(0) = u(0)$, $w(r) \geq u(r)$ and  
      $$F[w] (s)\leq -a (m^\prime)^{1+\alpha } {N-1\over s} \leq f(s)\ \mbox{on}\ [0, \delta].$$ 
Then by   Proposition \ref{comp1}, $w(s) \geq u(s)$ and then for all $s\in [0, \delta]$,  
$${u(s)-u(0)\over s} \leq -m^\prime.$$
This implies that $\Lambda_d(0) =\lambda_d(0)$. And the existence of the derivative at zero is proved.

We still need to prove that it is zero. Suppose by contradiction that it is not, one can suppose that it is positive, the other case can be done  with obvious changes. Let $\delta$ be small enough that 
$${a(u^\prime (0))^{1+\alpha} \over \delta 2^{1+\alpha} } > f(r) \ \mbox{for} \ r< \delta.$$ 
The function
 $$\varphi(x) = u(0) + \frac{u^\prime(0)}{2}x_1+{u^\prime (0)\over 4\delta}x_1^2.$$
touches $u$   from below  at zero and can be used as a test function and by definition of a  supersolution
  $${a (u^\prime (0))^{1+\alpha}\over \delta 2^{1+\alpha}}\leq F[\varphi](0)\leq f(0),$$
a contradiction with the choice of  $\delta$. This ends the proof.
   
\bigskip

We now state the main result of this section. 
\begin{theo}\label{regC1}
Suppose that $u$ is a radial solution of $F[u]= f$. Then  $u$ is  ${\cal C}^1$.  
  \end{theo}
The proof of  Theorem \ref{regC1} relies on Proposition \ref{propreg2}, Corollary \ref{deriv} and  Proposition \ref{propreg3}.

  \begin{prop}\label{propreg2}
Suppose that $u$ is a radial solution of $F[u] = f$
such that in $\bar r$ one of the derivative  numbers is different from zero. 

Then in  a 
neighborhood of $\bar r$, $u^\prime$ exists in the classical sense and the function $u^\prime$ is 
continuous in $\bar r$. 
\end{prop}
\begin{cor}\label{deriv}
If, at $r_1$, one of the derivative numbers is zero , then $u^\prime (r_1)$ exists and it is zero.
\end{cor}
{\em Proof of Proposition \ref{propreg2}.}
There are, in theory, 8 cases to treat, because each of the derivative number may be 
either positive or negative.
 But in fact  considering the function $v=-u$, that satisfies
$$G[v]:=-F[-v]=-f,$$
it is enough to consider only half of the cases.

What we want to prove is that,  for any $\bar r$,  as in Proposition \ref{propreg2}
\begin{equation}\label{thesis}
\exists\delta>0, \lambda_d(r)=\lambda_g(r)=\Lambda_d(r)=\Lambda_g(r),\  \forall r\in (\bar r-\delta, \bar r+\delta).
\end{equation}
{\bf Claim 1: $\Lambda_d(\bar r)=k<0$ implies the thesis (\ref{thesis}).}

We first prove that 
\begin{equation}\label{C11} 
 \Lambda_d(\bar r)<\mu<0\Rightarrow \exists\delta>0, \ \sup(\lambda_d,\lambda_g,\Lambda_d,\Lambda_g)(r) <\mu<0,\   \forall r\in (\bar r-\delta, \bar r+\delta).
\end{equation}
Let $\mu^\prime \in ]\Lambda_d, \mu[$, and  
  let $\delta_1$ be small enough in order that ${\mu^\prime \over 1-\sqrt{\delta_1}} > \Lambda_d$. Let $\delta_2 < \delta_1$ be small enough in order that 
 \begin{equation}\label{P1}{a |\mu|^{1+\alpha} \over (1-\sqrt{\delta_2})  \sqrt{\delta_2}}\geq |f|_\infty.
\end{equation}
Let $\delta_3 < {\bar r\over 2}$ such that for $\delta <\inf( \delta_1, \delta_2, \delta_3)$
   $${u(\bar r+\delta)-u(\bar r)\over \delta} \leq {\mu^\prime\over 1-\sqrt{\delta_1}}.$$
Fixe such a  $\delta$,  by continuity of $u$,  there exists $\delta_4$ such that for 
$r\in ]\bar r-\delta_4, \bar r+ \delta_4[$
\begin{equation}
\label{zxc}{u( r+\delta)-u( r)\over \delta} \leq {\mu\over 1-\sqrt{\delta_1}}.\end{equation}
For any such $r$  let 
      $$w := w ( \mu   ,{ \mu \over  (1-\sqrt{\delta})  \sqrt{\delta}} ,  r).$$
Then 
      $$w(r) = u(r), \ w(r+\delta) \geq u(r+\delta), $$
and using (\ref{P1}) it is easy to check that $w$ is a supersolution in $[r, r+\delta]$. 
 From this, using Proposition \ref{comp1}, one gets   that $w(s) \geq u(s)$ on $[r, r+\delta]$ and  then  
      $${u(s)-u(r)\over s-r} \leq \mu$$
which implies that  $\Lambda_d (r) \leq \mu$.  Exchanging the roles of $r$ and $s$ one 
obtains also $\Lambda_g (r) \leq \mu$.  This proves (\ref{C11} ).

To complete the proof of {\bf Claim 1} we prove that 
\begin{equation}\label{C12}
0> \Lambda_d(\bar r)>\mu \Rightarrow \exists\delta>0,\ \inf(\lambda_d,\lambda_g,\Lambda_d,\Lambda_g)(r) >\mu, \  \forall r\in (\bar r-\delta, \bar r+\delta).
\end{equation}
Indeed (\ref{C11}) and (\ref{C12}) imply the thesis (\ref{thesis}).

The proof of (\ref{C12}) proceeds in a similar fashion then that of (\ref{C11}), we give the detail of the computation for completeness sake. Let $\delta_1< \inf({\bar r\over 2},1) $ be small enough in order that ${\mu\over 1+\sqrt{\delta_1}}< \Lambda_d$,   that
\begin{equation}\label{P2}
 {8A(N-1)\over \bar r}<{a\over\sqrt{\delta_1}}\ \mbox{and that}\ 
{|\mu|^{1+\alpha}   a \over 2 (1+\sqrt{\delta_1})^{1+\alpha} \sqrt{\delta_1}}> |f|_\infty.
\end{equation}
 As above there exists $\delta_4$ such that  for $\delta< \delta_1$ small enough and 
 for $r\in ]\bar r-\delta_4, \bar r+\delta_4[$,
             $${u(r+\delta)-u(r)\over \delta} > {\mu\over 1+\sqrt{\delta_1}}.$$
We  define, on $[r, r+\delta]$,
$$w:=  w( \mu,    {|\mu|\over (1+\sqrt{\delta}) \sqrt{\delta}} , r)$$
then $w(r)= u(r), \ w(r+\delta) \leq u(r+\delta)$ and, using (\ref{P2}), $w$ is a sub-solution. Then                    
$$w(s) \leq u(s)\ \mbox{on}\ [r, r+\delta].$$
Finally one gets that  there exists $\delta_4$, 
and $\delta$ such that for $r\in ]\bar r-\delta_4, \bar r+\delta_4[$,  and $s\in [r, r+\delta]$
$$u(s) \geq u(r) + \mu (s-r),$$
hence  
$$\lambda_d (r)\geq \mu.$$
Exchanging the roles of $r$ and $s$, one gets $\lambda_g (r) \geq \mu$. This proves 
(\ref{C12}) and it ends the proof of {\bf Claim 1}.

Observe that in fact, (\ref{C11}) and (\ref{C12}) prove more than  {\bf Claim 1} because they imply that the derivative of $u$ is continuous in points where the derivative is not zero.

\medskip
{\bf Claim 1} implies 

\noindent{\bf Claim 2: $\lambda_d(\bar r)=k>0$ implies the thesis (\ref{thesis}).}

Indeed for $v=-u$, 
$$\lambda_{d,u}(\bar r)=-\Lambda_{d,v}(\bar r)$$ 
(here we have added in the index the function for which the derivative number is computed). 

\medskip
The proofs of the following claims are similar to the proof of  (\ref{C11}) or (\ref{C12}) but, of course, each case needs a different choice of function $w$, so we give the details for completeness sake.

\noindent{\bf Claim 3: $\Lambda_d(\bar r)=k>0$ implies the thesis (\ref{thesis}).}
 
                         Suppose that $0<\mu<\Lambda_d$,  and let $\delta_1$ be so that ${\mu^{1+\alpha}   a\over  (1-\sqrt{\delta_1})\sqrt{ \delta_1}} \geq |f|_\infty$. Let $\delta_2< \delta_1$ be so that ${\mu\over 1-\sqrt{\delta_2}} < \Lambda_d$. As it is done before  let $\delta_3<\inf (\delta_1, \delta_2)$ be so that for $\delta < \delta_3$, 
                         $${u(\bar r+\delta)-u(\bar r)\over \delta}  > {\mu\over 1-\sqrt{\delta_2}}.$$
Finally let  $\delta_4$ so that for $r\in [\bar r-\delta_4, \bar r+\delta_4]$ one has 
${u(r+\delta)-u(r)\over \delta} > {\mu\over 1-\sqrt{\delta_2}}$. Then define 
$w:=w(\mu,  {\mu\over (1-\sqrt{\delta} )\sqrt{\delta}}, r)$. Then 
$w(r) = u(r)$, $w(r+\delta) \leq u(r+\delta)$ and $w$ is a sub-solution. 

Hence,  by Proposition \ref{comp1},
$w(s)\leq u(s)$ on $[r,r+\delta]$, which  implies that  for 
$r\in ]\bar r-\delta_4, \bar r+\delta_4[$, 
$$\inf (\Lambda_d, \Lambda_g, \lambda_g, \lambda_d)(r) \geq \mu.$$
We are back to the hypothesis of {\bf Claim 2}, which implies the thesis .

\noindent{\bf Claim 4: $\Lambda_g(\bar r)=k<0$ implies the thesis (\ref{thesis}).}

Let $\mu$ be such that $0>\mu>\Lambda_g(\bar r)$.  Let $\delta_1 < {\bar r\over 2}$, 
$\sqrt{\delta_1} < \inf (1, {a\bar r\over 8A (N-1)})$, $ \Lambda_g(\bar r)<\frac{\mu}{1-\sqrt{\delta_1}}$, 
 and such that 
 $${|\mu|^{1+\alpha} a\over 2 (1-\sqrt{\delta_1})\sqrt{\delta_1}} > |f|_\infty.$$ 
  As before there exists $\delta_4<\delta_1$  such that for $r\in ]\bar r-\delta_4, \bar r+\delta_4[$  and for $\delta < \delta_4$ one has 
  $${u(r-\delta)-u(r)\over -\delta} \leq {\mu\over 1-\sqrt{\delta_1}}.$$
Then 
$$w(s) := w ( \mu   ,{ |\mu| \over  (1-\sqrt{\delta})  \sqrt{\delta}} ,  r)(s)$$
is  a subsolution which 
satisfies $w(r)=u(r)$, $w(r-\delta)\leq u(r-\delta)$. 
Then, by Proposition \ref{comp1},  $w(s)\leq u(s)$ in $(r-\delta,r)$.
This in turn  implies  that
$$ \sup(\lambda_d,\lambda_g,\Lambda_d,\Lambda_g)(r) <\mu<0.$$
We are once again in the hypothesis of {\bf Claim 1} and we are done.

\noindent{\bf Claim 5: $\Lambda_g(\bar r)=k>0$ implies the thesis (\ref{thesis}).}

Let $\mu$ be so that $\Lambda_g > \mu >0$. Let $\delta_1$  be such that  $\delta_1< {\bar r\over 2}$, $ {\mu\over 1-\sqrt{\delta_1}} <\Lambda_g$, 
$$\sqrt{\delta_1}  < \inf (1, {a\bar r\over 8A(N-1) } )\ \mbox{ and } {\mu^{1+\alpha}   a\over 2(1-\sqrt{\delta_1})\sqrt{ \delta_1}} \geq |f|_\infty.$$ 
As before for $\delta$ fixed, $\delta < \delta_1$,  there exists some $\delta_4< \delta_1$ 
such that for $r\in [\bar r-\delta_4, \bar r+\delta_4]$, 
$${u(r-\delta)-u(r)\over -\delta} \geq {\mu\over 1-\sqrt{\delta_1}}.$$
We define $ w := w(\mu, {-\mu\over(1-\sqrt{\delta})\sqrt{ \delta}}, r)$, then 
$w(r)= u(r), \ w(r-\delta) \geq u(r-\delta)$ and  by the assumptions 
$w$ is a supersolution  on $(r-\delta, r)$.   
Once more this implies that
$$\inf(\lambda_d,\lambda_g,\Lambda_d,\Lambda_g)(r) >\mu>0,\  \forall r\in (\bar r-\delta, \bar r+\delta).$$ 
And we conclude with {\bf Claim 2}.

\medskip

\noindent Again, using $v=-u$, the  {\bf Claims 3, 4 and 5} give that respectively
$\lambda_d(\bar r)=k<0$ or  $\lambda_g(\bar r)=k<0$ or  $\lambda_g(\bar r)=k>0$ imply 
the thesis (\ref{thesis}). And this ends the proof.

\bigskip

\noindent {\em Proof of Corollary \ref{deriv}.} By Proposition \ref{propreg2} if one of the derivative number has a sign then $u^\prime(r_1)$ exists and it is different from zero, which contradicts our hypothesis, so all four  derivative numbers are zero and $u^\prime(r_1)=0$ in the classical sense.

\bigskip
We finally give the last step of the proof of Theorem \ref{regC1}:
\begin{prop}\label{propreg3}
      $u^\prime$ is continuous on points where $u^\prime =0$. 
\end{prop}
{\em Proof.} We treat separately the case $\bar r=0$ and $\bar r\neq 0$. 
In the latter case,  let $\epsilon>0$  and $\delta_1<\inf ({\bar r\over 2},{a\epsilon^{1+\alpha}\over 2^{1+\alpha} |f|_\infty}) $, 
such that for $\delta< \delta_1$, ${u(\bar r-\delta)-u(\bar r)\over (-\delta)} \leq {\epsilon\over 4}$.

Fixing such $\delta$,  
let $\delta_2< \delta_1$ such that for $r\in [\bar r-\delta_2, \bar r+\delta_2]$,  
by continuity,   ${u(r-\delta)-u(r)\over (-\delta)} \leq {\epsilon\over 2} $.  Then 
$$u(r-\delta)\geq u(r) +{ \epsilon\over 2}  (-\delta).$$
Let $w:= w( \epsilon, {\epsilon\over  2 \delta}, r)$. 
Then $w(r)= u(r),\ w(r-\delta) \leq u(r-\delta)$ and, by the hypothesis on $\delta$, $w$ is a 
sub-solution on $[r-\delta, r]$. Then $w(s) \leq u(s)$ for  $r\in [\bar r-\delta_2, \bar r+\delta_2]$,
and $s\in [r-\delta, r]$. By passing to the limit when $s$ goes to $r$ it gives
$$u^\prime(r) \leq \epsilon.$$

We now prove that   for all $\epsilon>0$, there exists a  neighborhood of $\bar r$ where $u^\prime \geq -\epsilon$.

 Let, as above,  
 $$\delta_1<\inf ({\bar r\over 2},{a\epsilon^{1+\alpha}\over 2^{1+\alpha}  |f|_\infty}),$$ 
 such that for $\delta< \delta_1$, ${-\epsilon\over 4}\leq {u(\bar r-\delta)-u(\bar r)\over (-\delta)}$.

 Fixing such $\delta$,  let $\delta_2< \delta_1$ such that for $r\in [\bar r-\delta_2, \bar r+\delta_2]$,  by continuity,  
                  $${u(r-\delta)-u(r)\over (-\delta)} \geq {-\epsilon\over 2}.$$
Let $w:= w (-\epsilon, {-\epsilon\over  2\delta}, r)$. 
Then $w(r-\delta) \geq u(r-\delta)$, $w(r) = u(r)$ and $w$ is a super-solution on $[r-\delta, r]$. Using again Proposition \ref{comp1} 
$$w(s) \geq u(s), $$
and passing to the limit $u^\prime (r) \geq -\epsilon.$

\bigskip
\noindent We now consider the case where $\bar r = 0$. We want to prove the inequality $|u^\prime| \leq \epsilon$ in a neighborhood of zero. 
           
Take any $\epsilon>0$, by Lemma \ref{lemuprime0} there exists $\delta_\epsilon>0$ such that
$$ \left |{u(\delta)-u(0)\over \delta}\right| \leq {\epsilon\over 2},\ \forall \delta\in (0,\delta_\epsilon).$$
Let $\delta\leq \min( {a\epsilon^{1+\alpha} (N-1)\over 2|f|_\infty},\delta_\epsilon)$,
 by the continuity of $u$, there exists $\delta_1< \delta$ such that for $r\in [\delta, \delta+\delta_1]$, 
$${u(r)-u(r-\delta)\over \delta} \geq -\epsilon.$$
             
Let us consider   $w:= w( -\epsilon ,0,r)$.  Then 
$u(r)= w(r), \ w(r-\delta) \geq u(r-\delta)$ and, with our choice of $\delta$, $w$ is a supersolution  in $[r-\delta ,r ]$.  
By Proposition \ref{comp1} we have obtained that 
               $$u(s)\leq w(s) \ \mbox{on}\ [r-\delta ,r ]$$
and then
                 $$u^\prime(r) \geq -\epsilon, \ \forall r\in  [\delta, \delta+\delta_1].$$
In particular, $u^\prime (\delta)  \geq -\epsilon\ $ for any
$\delta\leq \min( {a\epsilon^{1+\alpha} (N-1)\over 2|f|_\infty},\delta_\epsilon)$.  
\medskip            

In a similar fashion, for any
$\delta\leq \min( {a\epsilon^{1+\alpha} (N-1)\over 2|f|_\infty},\delta_\epsilon)$ and for all $r\in [\delta, \delta+ \delta_1]$ 
  $$u^\prime (r) \leq \epsilon.$$
This is the 
desired result and it ends both the proof of Proposition \ref{propreg3} and the proof of Theorem \ref{regC1}. 

\section{${\cal C}^{1,\beta}$ regularity. }
 Observe that we now know that $u^\prime$ is continuous, so we can consider 
$\tilde F(x, D^2 v) := F(x, \nabla u, D^2 v)$
and clearly $u$ is a solution of
$$\tilde F(x, D^2 v):=f(x).$$

\begin{theo}\label{regC2}
Suppose that $N\leq 3$ or for any $N>3$ that $M\mapsto F(x, p, .)$ is convex or concave
and  that $u$ is a radial solution of $F[u]= f$. Then  $u$ is ${\cal C}^{1,{1\over 1+\alpha} }$
everywhere and is ${\cal C}^2$ on points where the derivative is different from zero. 
\end{theo}
The case where the derivative is different from zero is easy to treat:
\begin{prop}\label{C1beta}
 Suppose that for $r_o>0$, $u^\prime (r_o) \neq 0$ in the viscosity sense.  
Then,    on a neighborhood around $r_o$, $u$ is ${\cal C}^{1, \beta}$  for some $\beta$, and if  $N\leq 3$  or for any dimension when  $\tilde F$ is convex or concave,  $u$ is ${\cal C}^{2,\beta}$
 \end{prop}
{\em Proof  of Proposition \ref{C1beta}.}
Observe that condition [(H4)] implies that we are in the hypothesis  of \cite{Ca} in  $B_{r_o+\delta}\setminus \overline{B}_{r_o-\delta}$ for some $\delta>0$. Hence $u$ is ${\cal C}^{1,\beta}$ on  that  annulus.
Similarly if $F$ is concave or convex we are in the hypothesis of \cite{E} and \cite{Ca}  and
$u$ is ${\cal C}^{2,\beta}$.

Note that when $N\leq 3$,
the ${\cal C}^{2,\beta}$ regularity holds without the convexity  or concavity assumption thanks to \cite{NV}.

\bigskip
To  prove the ${\cal C}^{1, {1\over 1+\alpha}} $ regularity result on any point, including those with the derivative equal to zero,  we begin to establish a
technical lemma. 
\begin{lemme}\label{lem22}
Suppose that $u$ is a  radial ${\cal C}^2$ viscosity solution of (\ref{Theq}), and that on 
$]r, s[$, $0< r < s$, $u^\prime >0$.  Then 
 \begin{equation}\label{eqA}
        |u^\prime |^\alpha u^\prime( s) \leq|u^\prime |^\alpha u^\prime (r)  +  {(1+\alpha)} \int_r^s  \epsilon_{a,A}(f(t) ) dt
             \end{equation}
where $\epsilon_{a,A} (x) = {x^+\over a}-{x^-\over A}$;  furthermore, for $\gamma={A\over a}(N-1)(1+\alpha)$
              \begin{equation}\label{eqB}
        |u^\prime |^\alpha u^\prime( s) \geq \left({r\over s}\right)^{\gamma}|u^\prime |^\alpha u^\prime (r)  -{|f|_\infty   (1+\alpha)s\over A(N-1) (1+\alpha)+ A}  \left(1-\left({r\over s}\right)^{\gamma+1}\right).
             \end{equation}
If $u^\prime <0$  on $]r, s[$, $0< r < s$,  then  
\begin{equation}\label{eqC}
        |u^\prime |^\alpha u^\prime( s) \geq |u^\prime |^\alpha u^\prime (r)  +  {(1+\alpha)} \int_r^s  \epsilon_{A,a}(f(t) )dt
\end{equation} and 
\begin{equation}\label{eqD}
       |u^\prime |^\alpha u^\prime( s) \leq \left({r\over s}\right)^{\gamma}|u^\prime |^\alpha u^\prime (r)  +  {|f|_\infty(1+\alpha)s\over A(N-1)(1+\alpha)+a} \left(1- \left({r\over s}\right)^{\gamma+1 }\right).
             \end{equation}
\end{lemme}
{\em Proof }:
We start by supposing that $u^\prime>0$ in $(r,s)$.
Since $u^{\prime\prime}$ is continuous then $(u^{\prime \prime})^{-1}(\R^+)$ is  an open set of $\R^+$ and there exists  a union of numerable open sets $\cup_{n\in \N} ]r_n,r_{n+1}[$, with $]r,s[ = \cup_{n\in \N} ]r_n,r_{n+1}[$ where $u^{\prime\prime}$ is of constant sign on each interval $]r_n, r_{n+1}[$. 
By redefining in an obvious fashion the end points of the intervals  
one can suppose that $]r, s[ =  \cup_{n\in \Z}]r_n, r_{n+1}]$ with  
$u^{\prime \prime} \geq 0$ on $[r_{2p}, r_{2p+1}]$ and  $u^{\prime \prime} \leq 0$ on 
$[r_{2p+1}, r_{2p+2}]$. 

 Then, in     $[r_{2p}, r_{2p+1}]$,

\begin{equation}\label{asa}a(u^{\prime\prime} + {(N-1)\over r} u^\prime ) |u^\prime |^\alpha \leq f(r)\leq A(u^{\prime\prime} + { (N-1)\over r} u^\prime ) |u^\prime |^\alpha\ \end{equation}
and, in     $[r_{2p+1}, r_{2p+2}]$,
\begin{equation}\label{asb}(A u^{\prime\prime} + {a (N-1)\over r} u^\prime ) |u^\prime |^\alpha \leq f(r)\leq(a u^{\prime\prime} + {A (N-1)\over r} u^\prime ) |u^\prime |^\alpha .
\end{equation}
 We begin to prove (\ref{eqA}) using the inequality on the left of $f$.   
In each case one can drop the term $|u^\prime |^\alpha u^\prime$, hence integrating and 
using $\displaystyle{f\over A} , {f\over a} \leq \epsilon_{a,A} (f)$, this imply (\ref{eqA}) on 
$[r_{2p}, r_{2p+1}]$ and  on $[r_{2p+1}, r_{2p+2}]$. 

Let now $P$ arbitrary large negative  and $N$ arbitrary large positive,  $r_{P}$ close to $r$ and $r_N$ close to $s$, 

                 \begin{eqnarray*}
                  |u^\prime |^\alpha u^\prime(r_N) &\leq&|u^\prime |^\alpha u^\prime (r_{N-1})  +  {(1+\alpha)} \int_{r_{N-1}}^{r_N}  \epsilon_{a,A}(f(t))  dt\\
                  &\leq &  |u^\prime |^\alpha u^\prime (r_{P} )+ (1+\alpha)\sum_{P}^N \int_{r_{n}}^{r_{n+1}}   \epsilon_{a,A}(f(t) ) dt
                  \end{eqnarray*}
and one obtains (\ref{eqA}) by passing to the limit when $P$ and $N$ go 
respectively  to $-\infty$ and $+\infty$. 

We now prove (\ref{eqB}).    
In $[r_{2p+1}, r_{2p+2}]$, since $u^{\prime\prime}\leq 0$ we multiply the second inequality  of (\ref{asb}) by  $\frac{(1+\alpha)}{a}r^{A(N-1)(1+\alpha)\over a}:=\frac{(1+\alpha)}{a}r^\gamma$ and integrating one gets 
$$ r_{2p+2}^{\gamma} |u^\prime |^\alpha u^\prime (r_{2p+2}) \geq r_{2p+1}^{\gamma} |u^\prime |^\alpha u^\prime (r_{2p+1}) + \int_{r_{2p+1}}^{r_{2p+2}} {(1+\alpha)\over a}f(t) t^{\gamma} dt.$$
Hence dividing by $r_{2p+2}^{\gamma}$, using $f\geq -|f|_\infty$ one gets 
\begin{eqnarray*}
|u^\prime |^\alpha u^\prime (r_{2p+2})
& \geq&  \left({r_{2p+1}\over r_{2p+2}}\right) ^{\gamma}|u^\prime |^\alpha u^\prime (r_{2p+1}) - {(1+\alpha) \over a}{|f|_\infty}\int_{r_{2p+1}}^{r_{2p+2}} \left({t\over r_{2p+2}}\right)^{\gamma} dt.
\end{eqnarray*}
Similarly, if  $u^{\prime \prime} >0$, multiplying (\ref{asa}) by $(1+\alpha)r^{(N-1)(1+\alpha)}:=(1+\alpha) r^{\gamma_1}$ one gets
\begin{eqnarray*}
|u^\prime |^\alpha u^\prime (r_{2p+1})
&\geq& \left({r_{2p}\over r_{2p+1}}\right) ^{\gamma}|u^\prime |^\alpha u^\prime (r) - {(1+\alpha) \over A}{|f|_\infty}\int_{r_{2p}}^{r_{2p+1}}  \left({t\over r_{2p+1}}\right)^{\gamma_1} dt\\
&\geq&
\left({r_{2p}\over r_{2p+1}}\right) ^{\gamma}|u^\prime |^\alpha u^\prime (r_{2p}) - {(1+\alpha) \over a}{|f|_\infty}\int_{r_{2p}}^{r_{2p+1}}  \left({t\over r_{2p+1}}\right)^{\gamma} dt. 
 \end{eqnarray*}
We have used the fact that $\gamma_1=(N-1)(1+\alpha)< (N-1)(1+\alpha ) {A\over a}:=\gamma$ and $A(N-1)(1+\alpha)+ A\geq A(N-1)(1+\alpha)+a$.

 Using the same decomposition of $]r,s[ = \cup_{n\in\Z} ]r_n, r_{n+1}]$, with $u^{\prime\prime}$ of constant sign in each interval,  for $P$ large negative and $N$ 
large positive , $r_{P}$ close to $r$ and $r_N$ close to $s$, one has 
\begin{eqnarray*}
                 |u^\prime |^\alpha u^\prime (r_N) &\geq&  \left({r_{N-1}\over r_{N}}\right) ^{\gamma}|u^\prime |^\alpha u^\prime (r_{N-1})  -{|f|_\infty  (1+\alpha) \over a} \int_ {r_{N-1}}^{r_N} \left({t\over r_N}\right) ^{\gamma} dt\\
            &\geq &
               \left({r_{N-2}\over r_{N-1}}\right) ^{\gamma} \left({r_{N-1}\over r_{N}}\right) ^{\gamma}|u^\prime |^\alpha u^\prime (r_{N-2}) \\
               &&-  \left({r_{N-1}\over r_N}\right) ^{\gamma}{|f|_\infty   (1+\alpha) \over a}
\int_ {r_{N-2}}^{r_{N-1}}                 \left({t\over r_{N-1}}\right) ^{\gamma} dt\\
&&-         {|f|_\infty (1+\alpha) \over a}     \int_ {r_{N-1}}^{r_N} \left({t\over r_N}\right) ^{\gamma} dt\\
 &=&     \left({r_{N-2}\over r_N}\right) ^{\gamma}|u^\prime |^\alpha u^\prime (r_{N-2}) -    {|f|_\infty   (1+\alpha) \over a}\int_ {r_{N-2}}^{r_N} \left({t\over r_N}\right) ^{\gamma} dt\\
 &\geq & \left({r_{P}\over r_N}\right)^{\gamma}|u^\prime |^\alpha u^\prime (r_{P}) - {|f|_\infty  (1+\alpha) \over a}  \int_ {r_{P}}^{r_N} \left({t\over r_N}\right) ^{\gamma} dt.
 \end{eqnarray*}
 By passing to the limit   when $P$ and $N$ go to $-\infty$ and $+\infty$ one obtains (\ref{eqB}) .

The  inequalities (\ref{eqC}) and (\ref{eqD}) can of course  be  proved either in the same
manner or considering $v=-u$ as the solution of 
                      $$G[v] = -f$$
and $G[v] = -F[-v]$ which possesses the same properties as $F$.

   \begin{prop}
    The solutions of $F[u]= f$ are ${\cal C}^{1,{1\over 1+\alpha}} $.
\end{prop} 
\noindent   {\em Proof.} 
Let $r_1 > 0$ such that $u^\prime(r_1)=0$, and let ${r_1\over 2} < r< r_1$.  We shall prove that 
 \begin{equation}\label{machin}
     |u^\prime|^{\alpha +1} (r) \leq {2^{\gamma-1}(\gamma+1) |f|_\infty(1+\alpha)\over A} (r_1-r).
 \end{equation}
For that aim, suppose first that $u^\prime (r) >0$ and let $s$ be the first point 
between $r_1$ and $r$, so that $u^\prime (s) = 0$. Then $u^\prime >0$ between 
$s$ and $r$ and inequality (\ref{eqB}), with $\gamma= {A(N-1)(1+\alpha)\over a}$, becomes
$$ |u^\prime |^\alpha u^\prime( s)=
0 \geq \left({r\over s}\right)^{\gamma}|u^\prime |^\alpha u^\prime (r)  -{|f|_\infty   (1+\alpha)s\over A(N-1) (1+\alpha)+ a}  \left(1-\left({r\over s}\right)^{{\gamma}+1}\right) $$
and then 
\begin{eqnarray*}
|u^\prime |^\alpha u^\prime (r)&\leq &\left({s\over r}\right)^{\gamma} {|f|_\infty   (1+\alpha)\over A(N-1) (1+\alpha)+ a}  \left(1-\left({r\over s}\right)^{{\gamma}+1}\right)\\
                 & \leq&  ({\gamma}+1){s^{\gamma-1} \over r^{\gamma-1} }{|f|_\infty   (1+\alpha)\over A(N-1) (1+\alpha)+ a}  (s-r)\\
                 &\leq & 2^{\gamma-1} ({\gamma}+1) {|f|_\infty   (1+\alpha)\over A(N-1) (1+\alpha)+ a}  (r_1-r).
\end{eqnarray*}
 From this one gets that 
                  $$|u^\prime(r) |\leq C (r_1-r)^{1\over 1+\alpha}.$$
The case where $u^\prime (r)<0$ can be done similarly. 
                  
We now consider the right of $r_1$.  This proof still holds when $r_1=0$. Suppose $s> r_1$,  we want to prove that 
$$|u^\prime|^{\alpha+1} (s) \leq {(1+\alpha) |f|_\infty\over a}   (s-r).$$ 
 For that aim suppose that $u^\prime (s) >0$. Let $r$ be the last point  in $]r_1, s[$ such that $u^\prime (r) = 0$. Then 
 $u^\prime>0$ on $]r, s[$ and by (\ref{asa})
 $$ |u^\prime|^\alpha u^\prime (s) \leq 0+ (1+\alpha) \int_r^s \epsilon_{a, A} (f) .  $$
  This implies 
  
  $$|u^\prime|^{\alpha+1} (s) \leq  {(1+\alpha) |f|_\infty\over a}   (s-r)\leq  {(1+\alpha) |f|_\infty\over a}   (s-r_1). $$

\section{Appendix : The equivalence of definitions of viscosity solution in the case $\alpha >0$.}
We have the following 
equivalence result 
         
 \begin{prop}\label{propalphapo} If $F$ satisfies (H1) and (H2) with $\alpha\geq 0$,
the viscosity solutions in the classical meaning are the same   as 
the viscosity solutions in the sense of Definition \ref{defvs}.
\end{prop}
           
\begin{rema}
Let us note that  the Definition \ref{defvs} presents an advantage with regards to the 
classical definition since it allows to not test  points where the  gradient of a test function is 
zero.
\end{rema}
{\em Proof  of Proposition \ref{propalphapo}.}
              
We assume that $u$ is a supersolution in the sense of Definition \ref{defvs} and we 
want to prove that it is a supersolution  in the classical sense. We suppose that for $x_o\in\Omega$ there exists $M\in S$ such that 
\begin{equation}\label{mnb}
u(x)\geq u(x_o)+\frac{1}{2}\langle M(x-x_o),(x-x_o)\rangle +o(|x-x_o|^2):=\phi(x).\end{equation} 
Let us observe first that one can suppose that $M$ is invertible, since if it is not, it can be replaced by  $M_n= M-{1\over n} I $ which satisfies (\ref{mnb}) and tends to $M$.
Without loss of generality we will suppose that $x_o=0$.

Let $k> 2$ and $R>0$ such that
          $$\inf_{ x\in B(0, R)} \left({u(x)-{1\over 2} \langle Mx, x\rangle+ |x|^k } \right)= u(0)$$
where the infimum is strict. 
We choose $\delta< R$ such that $(2\delta)^{k-2} < {\inf _i |\lambda_i (M)|\over 2k }$. 
Let $\epsilon $ be such that 
$$\inf_{ |x|> \delta }  \left(u(x)-{1\over 2} \langle Mx, x\rangle+  |x|^k \right)= u(0)+ \epsilon $$
and let $\delta_2< \delta $ and such that $k(2\delta)^{k-1}  \delta_2 + |M|_\infty (\delta_2^2 + 2 \delta_2 \delta) < {\epsilon\over 4}$. 
Then, for $x$ such that $|x|< \delta_2$,  
           \begin{eqnarray*}
           \inf_{ |y|\leq \delta} \{ u(y)-{1\over 2} \langle M(y-x), y-x\rangle +  |y-x|^k\} &\leq&      \inf_{|y|\leq \delta} \{ u(y)-{1\over 2} \langle My, y\rangle +  |y|^k\}+ {\epsilon\over 4}\\
           &=& u(0) + \epsilon/4
           \end{eqnarray*}
and on the opposite 
\begin{eqnarray*}
\inf_{ |y|> \delta} \{ u(y)-{1\over 2} \langle M(y-x), y-x\rangle +  |y-x|^k\} 
&\geq&      \inf_{|y|> \delta} \{ u(y)-{1\over 2} \langle My, y\rangle +  |y|^k\}- {\epsilon\over 4}\\
&> &u(0) + {3\epsilon\over 4}.
\end{eqnarray*}
Since the function $u$ is supposed to be non locally constant, there exist $x_\delta$ 
and $y_\delta$ in $B(0, \delta_2)$ such that 
        $$ u(x_\delta) > u(y_\delta) -{1\over 2} \langle M (x_\delta-y_\delta), x_\delta-y_\delta\rangle +  |x_\delta-y_\delta|^k$$
         and then the infimum $\inf _{y, |y| \leq \delta}\{u(y) -{1\over 2} \langle M (x_\delta-y), x_\delta-y\rangle + 
  |x_\delta-y|^k\}$ is achieved on some point $z_\delta$ different from $x_ \delta$. This implies that the function 
$$\varphi(z):=u(z_\delta ) + {1\over 2} \langle M(x_\delta- z), x_\delta-z) - |x_\delta-z|^k- {1\over 2} \langle M(x_\delta-z_\delta), x_\delta-z_\delta \rangle +  |x_\delta -z_\delta|^k$$ 
touches $u$ by below at the  point $z_\delta$. 
But 
$$\grad\varphi(z_\delta)=M(z_\delta -x_\delta ) -k  |x_\delta-z_\delta|^{k-2} (z_\delta -x_\delta),$$    
cannot be zero since, if it was,  $z_\delta -x_\delta$ would be an eigenvector for the eigenvalue 
$k |x_\delta-z_\delta |^{k-2} $ which is supposed   to be strictly less than any eigenvalue of $M$. 

We have obtained that
$\grad\varphi(z_\delta)\neq 0$
and then, since $u$ is a supersolution in the sense of Definition \ref{defvs},
$$F( z_\delta, M(z_\delta -x_\delta) -k |x_\delta -z_\delta |^{k-2} (z_\delta -x_\delta ), M - {d^2\over dz^2} (|x_\delta-z|^k)(z_\delta) ) \leq g(z_\delta, u(z_\delta)).$$
By passing to the limit we obtain 
               $$0\leq g(0, u(0)),$$
 which is the desired conclusion.

Of course we can do the same for sub-solutions.


\begin{thebibliography}{99}

\bibitem{BNV} H. Berestycki, L.  Nirenberg, S.R.S. Varadhan, \emph{  The 
principal eigenvalue and maximum principle for second-order elliptic operators 
in general domains,} Comm. Pure Appl. Math. \textbf{47}  (1994),  no. 1, 47--92.

\bibitem{BD1} { I. Birindelli, F. Demengel,} 
{ \it Eigenvalue, maximum principle and regularity for fully non linear
homogeneous operators,} Comm. Pure and Appl. Analysis, {\bf 6} (2007), 335-366.



\bibitem{BD2}  I. Birindelli, F. Demengel, 
 {\em Regularity and uniqueness of the first eigenfunction for
singular fully nonlinear operators}, J. Differential Equations {\bf 249} (2010) 1089-1110. 
 
 \bibitem{BDr}  I. Birindelli, F. Demengel,  {\em Uniqueness of the first  eigenfunction for fully nonlinear equations: the radial case}, Journal for Analysis and its Applications, {\bf 29} (2010), 75-88.


\bibitem{CaC1} X. Cabr\'e, L. Caffarelli, {\em  Regularity for viscosity solutions of 
fully nonlinearequations $F(D^2u) = 0$,}  Topological Meth. Nonlinear Anal. {\bf 6} (1995), 31-48. 

\bibitem{CaC2} X. Cabr\'e, L. Caffarelli, {\em Interior ${\cal C}^2$ regularity theory for 
a class of nonconvex fully nonlinear elliptic equations} J.   Maths Pures Appl.  (9) {\bf 82} (2003) 573-612.

\bibitem{Ca} L. Caffarelli, {\em Interior a Priori Estimates for Solutions of Fully 
Nonlinear Equations}, 
The Annals of Mathematics, Second Series,  {\bf 130} (1989), 189-213.

\bibitem{CC} { L. Caffarelli, X. Cabr\'e,} Fully-nonlinear equations
Colloquium Publications 43, American Mathematical Society, Providence, RI,1995.

\bibitem {CDV} I. Capuzzo Dolcetta, A. Vitolo, 
{\em GlaeserÕs type gradient estimates for non-negative 
solutions of fully nonlinear elliptic equations } 
 Discrete Contin. Dyn. Syst. 28 (2010), 539-557.
 
\bibitem{DFQ}{G. Davila, P. Felmer, A. Quaas,} {\em  Harnack Inequality for singular fully 
nonlinear operators and some existence's results} .
Calculus of Variations and PDE, Vol. 39,(2010), 557-578.

\bibitem{E} L.C. Evans, {\it Classical Solutions of Fully Nonlinear, Convex,
Second-Order Elliptic Equations}, Comm. on Pure and Applied Mathematics, {\bf 25}, (1982)333-363 .

\bibitem{ES} L.C. Evans, O. Savin, {\em $C^{1,\alpha}$ Regularity for infinity 
harmonic functions in two dimensions.}   Calculus of Variations and PDE, {\bf 32} (2008),  325-347.
 
\bibitem{NV} N. Nadirashvili, S. Vladut, {\em On Axially Symmetric Solutions of Fully Nonlinear Elliptic Equations},  Math. Z. to appear.

\end{thebibliography}
       \end{document}